\documentclass[12pt]{amsart}

\usepackage{fullpage}
\usepackage{amsfonts}
\usepackage{enumerate}
\usepackage{verbatim}
\usepackage{graphicx}
\usepackage{bbold}

\title{The smoothness test for a density function}

\author{ Bogdan \'Cmiel }

\author{ Karol Dziedziul }
\author{ Barbara Wolnik }
\address{AGH University of Science and Technology}
\address{Gda\'nsk University of Technology}
\address{University of Gda\'nsk}

\newtheorem{corollary}{Corollary}[section]
\newtheorem{lemma}{Lemma}[section]
\newtheorem{thm}{Theorem}[section]
\newtheorem{Proposition}{Proposition}[section]
\newtheorem{definition}{Definition}[section]
\newtheorem{assumption}{Assumption}[section]
\newtheorem{remark}{Remark}[section]
\newtheorem{example}{Example}[section]
\newtheorem{ev}{Condition}

\newcommand{\R}{\mathbb{R}}
\newcommand{\supp}{\operatorname{supp}}

\newcommand{\N}{\mathbb{N}}
\newcommand{\Z}{\mathbb{Z}}

\newcommand{\E}{{\mathbf {E}}}

\begin{document}

\begin{abstract}
The problem of testing hypothesis that a density function has no more than $\mu$ derivatives versus it has more than $\mu$ derivatives is considered. For a solution, the $L^2$ norms
of wavelet orthogonal projections on some orthogonal "differences" of spaces from a multiresolution analysis is used. 
For the construction of the smoothness test an asymptotic distribution of a smoothness estimator is used. To analyze that asymptotic distribution, a new technique of enrichment procedure is proposed. The finite sample behaviour of the smoothness test is demonstrated in a numerical experiment in case of determination if a density function is continues or discontinues.
\end{abstract}


\maketitle

\hspace{0.5cm}

\section{Introduction}

The smoothness estimation problem was recently analyzed in \cite{Dziedziul} and \cite{CD}. In the first paper the smoothness estimator was obtained by histogram approach. The drawback of that method was the range of applications. It has to be assumed that a smoothness parameter of a density function is smaller than the smoothness of a piecewise constant function. In the second paper that uncomfortable assumption was removed using wavelet approach. In that case the  smoothness has to be only smaller than the smoothness of a compactly supported wavelet which was used in the estimation. That method allows us to concern all finite ranges of the smoothness parameter.
It was also shown that the estimator is "pseudo-consistent" on Besov classes of density functions. Pseudo-consistency means that the lower limit of an estimator, when the experiment size goes to infinity, is equal to the estimated
parameter with probability one (there is "$\liminf$" instead of "$\lim$"). In \cite{CD} it was also shown that if we restrict the Besov class to the "piecewise-smooth" class of functions then the smoothness estimator is strongly consistent.
\\
\\
In this paper we focus on the test of the hypothesis that the smoothness parameter of a density function is smaller than or equal to some real value against the hypothesis that it is greater than that value. To obtain the form of that test, at some significance level, we analyze the asymptotic distribution of the smoothness estimator using Berry-Esseen inequality. It turns out that in the proof of the Berry-Esseen inequality (for the smoothness parameter estimator) we have to make some technical assumption on the density function. To avoid another restriction for the class of the piecewise-smooth density functions we propose the enrichment procedure. Namely to a raw sample of size $n$ from a density $f$ we propose to add a sample from an appropriate density $\xi$ of size
$
n_1=\pi n/(1-\pi),
$
where $0<\pi<1$.
Hence we obtain a sample from the density
$
f_\pi=(1-\pi)f+\pi \xi,
$
of the size $n/(1-\pi)$. Since the smoothness of the density $\xi$ is greater than the range of the smoothness for function $f$ then the density $f_\pi$ has the same smoothness parameter as function $f$ and satisfies our assumption. We demonstrate the behavior of the enrichment procedure as well as the smoothness test in a numerical experiment.
\\
\\
The class of piecewise-smooth functions, considered in this paper, is in the center of interest for problems of finding "change points" defined as smoothness defects. That problem was analyzed in \cite{H} for a regression model. This class is also important in many others applications, see \cite{LL}. In our approach it is possible to test the smoothness locally (on some intervals) since we use wavelet bases. It allows us to introduce an alternative method for estimating locations and sizes of smoothness defects on which we will focus in the future.
\\
\\
Over the last few years many papers have referred to smoothness identification. Smoothness tests were recently studied in \cite{Belitser} where an empirical Bayes approach is used to test the smoothness of a signal in a Gaussian model. The smoothness parameter was defined in terms of Sobolev spaces. Some methods of detection of a function from anisotropic Sobolev classes are also considered in \cite{Ingster}, where the authors use Fourier coefficients but those hypotheses have a different structure than ours.
The smoothness test also appears in the context of confidence sets for an unknown density in $B^s_{\infty,\infty}(I)$, where $I$ is an interval (see \cite{Hoffmann} and \cite{Gine}), in $B^s_{2,\infty}(I)$ (see \cite{BullNickl}), or in $B^s_{\infty,\infty}(M)$, where $M$ is a homogeneous manifold (see \cite{Gerard}).
The Besov spaces  $B^s_{p,\infty}$  gives us an opportunity to define a continuous scale of smoothness so recently we have proposed a direct way to estimate the smoothness parameter (see \cite{CD} and \cite{Dziedziul}). Let us compare our approach with earlier methods which use wavelets i.e. $r$-regular wavelet basis ($r$-RWB see Definition \ref{r-RWB}).
    In papers  \cite{Hoffmann} and \cite{Gine} the following assumption is essential:  the density $f\in B^s_{\infty,\infty}$, where $s<r$ for $r$-RWB satisfies the following condition:
\begin{ev}\label{BE}
There exists $b$ such that for all $j\geq 0$
\[
b 2^{-js} \leq \|P_j f -f\|_p,
\]

\end{ev}
\noindent where $P_j$  denotes orthogonal projection connected with $r$-RWB (see  (\ref{Pj})).
It is shown that the class of functions that does not satisfy Condition \ref{BE} is nowhere dense in $B^s_{\infty,\infty}$ (see Proposition 4, \cite{Gine}). It is worth mentioning that similar results  hold for all $B^s_{p,\infty}$, $s\in \R_+$, $p\geq 1$.
Using Condition \ref{BE} confidence sets are constructed as well as a test. Unfortunately a separation from the level zero condition for densities is needed (a class $\mathcal D$ (3.3) in \cite{Gine}). Note that separation from the level zero is essential in  the idea of enrichment procedure (see Example \ref{przyklad} and the explanation below).
In paper  \cite{Gerard}
 Condition \ref{BE} also appears. In that paper a class of functions on sphere which satisfies Condition \ref{BE} is given explicitly (see Proposition 6, \cite{Gerard}).
In our paper we define a  class of  piecewise-smooth functions (Definition \ref{P-S}) that satisfies a little stronger condition than Condition \ref{BE}, i.e.
\begin{equation}
\label{stronger}
b 2^{-js} \leq \|Q_j f\|_p,
\end{equation}
and $Q_j$ denotes orthogonal projection connected with $r$-RWB (see (\ref{Qj}) and  Definition \ref{P-S}).
Note that indeed the inequality \eqref{stronger} is stronger than Condition \ref{BE}:
\[
\|Q_j\|_p=\|P_{j+1}-P_j \|_p=\|P_{j+1}-P_{j+1}P_j \|_p\leq \|P_{j+1}\|_p \|P_j-I\|_p\leq C_p \|P_j-I\|_p,
\]
where $I$ is identity operator.
\\
\\
Since we use a zero oscillation condition (formula (\ref{moment})) we have \eqref{stronger} for
 piecewise-smooth functions
 $f\in C^k$,  $k< d(r)$. This result is  stated in \cite{CD} for $k<r$ but a proof shows that it is true for $k<d(r)$.
In applications such property might be important since $r \approx 0.2*d(r)$ for Daubechies wavelets.
In other papers, authors can characterize only functions with smoothness less than $r$.
  Asymptotic theorems in Section \ref{psf} give us full characterization of $\|Q_j f\|_p$, for all $1<p<\infty$, regardless of smoothness parameter of function $f$ (see Theorem \ref{p}).
 \\
\\
Our paper can be treated as 
a detailed analysis of some testing problem, that appears in the Bull's and Nickl's paper \cite{BullNickl}, in some special class of density functions. 
In that paper, to construct a test \eqref{TEST}, authors used an U-statistic very similar to (\ref{LNJ}). Let $X_1,X_2,...,X_n$ be i.i.d. with common probability density $f$ on $[0,1]$. let $\Sigma$ be
any subset of a fixed Sobolev ball $\Sigma(t, B)$ for some $t > 1/2$ and consider testing
 \[ 
f\in \Sigma \ \ \ \mathrm{against}  \ \ \ 
 f\in  \Sigma(t,B)\setminus  \Sigma, \quad \|f-\Sigma \|_2\geq \rho_n,
 \] 
where $\rho_n \geq 0$ is a sequence of nonnegative real numbers  
 (see \cite{Car}). We can take  $\Sigma=\Sigma(s,B)$ where $s>t$. The test statistics in paper \cite{BullNickl} is given by
\[
\Psi_n=1\{\inf_{g\in \Sigma}|T_n(g)|>\tau_n\},
\]
where $(\tau_n)$ is a sequence of thresholds and $T_n(g)$ is the U-statistics indexed by all functions from the set $\Sigma$.
Our test is based on a similar idea but we consider
a particular situation: the unknown density function $f$ is a piecewise-smooth function. Thus we managed to present our test in an explicit form as the  inequality (\ref{test2}) with the exact constants.
We estimate $m=m(f)$ such that $f\in C^{m-1}(\R)$ and $f \notin C^m(\R)$ using a smoothness parameter approach. We say that $m(f)$ is an index of function $f$ and we use the notation $id(f)=m(f)$. We were inspired by the paper of Horvath and Kokoszka \cite{H}. The index of function $f$ gives us an exact smoothness parameter $s^*_2(f)=id(f)+1/2$ (see Definition \ref{smooth}). Since we have an asymptotic distribution of the estimator of $s^*_2(f)$  we
have in fact a test in more classic form $H_0:s^*_2(f)=k/2 $ against $H_1: s^*_2(f)\neq k/2$, where $k\in \N$. We use and explore a concept of $id(f)$ to present our results in a more intuitive way.
In our paper we take
a class of piecewise-smooth functions ${\mathcal S}={\mathcal S}([-1,1],d(r),\Delta_1,\Delta_2, N_S)$ (see Definition \ref{class}). Then  classes
\[
{\mathcal S}_m:=\{f\in {\mathcal S}: id(f)=m, f\geq 0, \int_R f=1 \}, \quad  m=0,1,2,...
\]
  are  separated in $L^\infty$ norm (in fact derivatives are separated in $L^\infty$ norm)  and
 \[
{\mathcal S}_m \subset B^{s_1}_{2,\infty}\cap B^{s_2}_{\infty,\infty}
\]
 where $s_1\leq m+1/2$ and ( $s_2\leq m$ if $m\geq 1$ ) (see Theorem \ref{p} and Definition \ref{class}).
 Note that
 we can take the parameters of $\mathcal S$ (or in other words a separation criteria in $L^\infty$ norm), i.e. $\Delta_1,\Delta_2, N_S$ which depend on sample $n$, as in \cite{BullNickl} but we want to avoid laborious calculations. Our classes are not separated in $L^2$ norm as in \cite{BullNickl}.
Since we consider piecewise-smooth functions we have a control on the constants 
which also depend on wavelets (see $F_\psi$ Example 7.1), so it is important which basis we take. Moreover since our methods depend on $d(r)$ (the zero oscillation parameter of wavelet) we obtain an extra freedom. We find this important  because if we increase $r$ (smoothness of a wavelet) a wavelet support also increases, so from numerical point of view, to detect a smoothness defect we need a greater resolution level. Note that our approach covers also discontinuous functions (see Example 7.1) which makes an important difference from results in \cite{BullNickl}. Since our class contains discontinuous functions, then it is not a special case of classes considered in \cite{BullNickl}. Another difference is that our test is ready for implementation and we have checked it's finite sample behaviour.
\\
\\
The smoothness test in a regression problem is considered in  paper \cite{Car}.  It is a modified and more detailed version of the test from \cite{BullNickl}. In both papers the asymptotic distributions are not considered. In our paper we use an asymptotic distribution of the smoothness estimator. It is more statistical approach.
Bull's and Nickl's method is connected with large deviation approach (Bernstein's inequality) or "sharp rates" as in \cite{Car}.  
 \\
 \\
In our paper we state some auxiliary results in $L^p$ norm  but to construct the test we use $L^2$ norm only. This point have to be emphasized. Our test concerns the smoothness parameter $s^*_2(f)$ only, not $s^*_p(f)$ for all $p>1$.
We use $L^2$ norm to make a conclusion for index of function   $id(f)+1/2=s^*_2(f)$, see Remark 3.4.  We have
\[
id(f)=\lim_{j\to \infty} -\frac{\log_2 \|Q_j f\|_2}{j}-\frac{1}{2}.
\]
(see Theorem \ref{p}).
 Recall that Assumption \ref{BE} is essential. Moreover we have the following  result (a reformulation of Proposition \ref{col}):
if $f$ is a  piecewise smooth function with $id(f)\leq d(r)$ then there is $C$ (the constant $C$ is precise but not the best) depending on $f$ and $r$-RWB such that
\begin{equation}
\label{obciazenie}
\left| -\frac{\log_2 \|Q_j f\|_2^2}{2j}-\frac{1}{2}- id(f)\right|\leq \frac{C}{j}.
\end{equation}
Let be given a sample $X_1,\ldots,X_n$ from a density $f$. We consider two consistent estimators of $id(f)$:
\[
\widehat{id(f)}_{n}=-\frac{\log_2 E_{n,j(n)}}{2j(n)}-\frac{1}{2}
 \]
 (see (\ref{ENJ}) and Theorem \ref{TwENJ})  or
 \[
\widehat{id(f)}_{n}=-\frac{\log_2 L_{n,j(n)}}{2j(n)}-\frac{1}{2}
 \]
(see (\ref{LNJ}) and Corollary 4.1).
The estimator $E_{n,j}$ was examined in \cite{CD}. It is  easier to use it in numerical simulations.
In this paper we consider mainly $L_{n,j}$ since it is easier to
formulate a concentration theorem for it (Theorem \ref{thm55}) as well as a test below. Both estimators are closely related.
 We do not consider as usual a bias and a variance of $\widehat{id(f)}_{n}$. Since $L_{n,j(n)}$ is unbiased estimator of $\|Q_jf\|^2_2$ i.e. $ \E L_{n,j(n)}=\|Q_jf\|^2_2$ then the precision of the
 estimation is obtained by  \eqref{obciazenie}.
\\
\\
 It is important to have an analogue of \eqref{BE}
but in $f(x) dx$ measure
 namely
 \[
 (b 2^{-js})^p \leq \int_{\R} |Q_j f(x)|^p f(x) dx =:\delta_j(f).
 \]
 Unfortunately this might not be true for some piecewise-smooth function $f$. So not to lose generality we propose an enrichment procedure.
 In this way we obtain $f_\pi$ (see Theorem \ref{enrichment}) with the same smoothness as $f$, i.e. $id(f)=id(f_\pi)$.
 In this way we avoid the extra assumption. As numerical calculations show, the estimator is located better when we use the enrichment procedure. Now using reformulated Theorem \ref{thm55} we can state that
 asymptotically with an appropriate choice of $j=j(n)$
\begin{equation}
\label{CLT}
L_{n,j(n)}=\|Q_j f_\pi\|_2^2 + \frac{2\delta_j(f_\pi)}{\sqrt{n}} U,
\end{equation}
where $U\sim N(0,1)$.
We have more since we use Berry Esseen inequality. We decide to use it since in Central Limit Theorem we cannot control a convergence. This gives us
a test (below) to verify a hypothesis
\begin{equation}
\label{TEST}
H_0: id(f)\leq\mu_0\ \ \mathrm{against}\ \ H_1: id(f)\geq \mu_0 +1.
\end{equation}
 A naive rule for the above test $H_0: id(f)\leq\mu_0$ against $H_1$ is the following:
 an acceptance of $H_0$ if $\widehat{id(f)}_{n} \leq \mu_0+\frac{1}{2}$ and a rejection if $\widehat{id(f)}_{n} > \mu_0+\frac{1}{2}$,
but this procedure does not
include any significance level $\alpha$.
 \\
\\
 It follows from (\ref{test2}) that a rejection rule might be written as
\begin{equation}
\label{regula}
\widehat{id(f)}_{n}=-\frac{\log_2 L_{n,j}}{2j} - \frac{1}{2}\geq \mu_0 -\frac{\log_2\left(z_\alpha C_1 2^{j (\mu_0+1/2)}/\left(\mu_0! \sqrt{n}\right)+C_2/(\mu_0!)^2  \right)}{2j},
\end{equation}
or using the formula $n\asymp 2^{j(2d(r)+3)}$
\[
\widehat{id(f)}_{n}\geq \mu_0 -\frac{\log_2\left(z_\alpha C_1/\left(\mu_0!\ n^{\frac{d(r)+2-\mu_0}{2d(r)+3}}\right)+C_2/(\mu_0!)^2 \right)}{\log_2 n^{\frac{2}{2d(r)+3}}}
\]
\[
=\mu_0+D(n,\mu_0,d(r),\alpha).
\]
where $C_1$ and $C_2$ are precise constants dependent on function class $\mathcal S$, $r$-RWB and the enrichment procedure on level $\pi$ (see (\ref{test2})).
We can also determine precisely a power of the statistical test. The procedure to determine a formula on the power of test is standard. We should use both sides of inequality in Proposition \ref{col}.  Note only that since $d(r)\geq \mu_0$ then $D(n,\mu_0,d(r),\alpha)\rightarrow 0$ when $n\rightarrow \infty $ and because of \eqref{CLT}  the power of the test converges to $1$. In the last section we show how well the test works.
\\
\\
 From an application point of view our test works well if one try to detect a continuity or a discontinuity of a density function (see last section).
If one tries to detect a higher defect, the size of a sample increases rapidly not only in our approach but in all approaches from papers \cite{BullNickl}, \cite{Gerard}, \cite{Hoffmann}, \cite{Gine}. Note that the case of detection  discontinuity is not included by any of those papers.

\section{Smoothness parameter}

Let us recall the definition of Besov spaces in terms of an isotropic modulus of smoothness on $\R$. We denote the norm in $L^p(\R)$ by $\|\cdot\|_p$ and the space of all continuous and bounded functions by $C_b (\R)$.

\begin{definition}\label{defbesov}
For $t>0$, $k\geq 1$, $0<s<k$, $1\leq p < \infty$, $1\leq q < \infty$ let us denote
\[
\Delta^k_h f(x)=\sum_{j=0}^k (-1)^{j+k} {k\choose j} f(x+jh),\ \ \ \ \omega_{k,p}(f;t)=\sup_{0<h\leq t} \|\Delta^k_{h} f\|_p\ ,
\]
\begin{equation}\nonumber
\|f\|^{(s)}_{p,q}=\|f\|_p+\left(\int_0^{\infty} \left(\frac{\omega_{k,p}(f;t)}{t^s}\right)^q \frac{dt}{t}\right)^{1/q}.
\end{equation}
Besov spaces are defined as follows:
\[
B^s_{p,q}(\R)=\{f\in L^p(\R): \|f\|^{(s)}_{p,q}<\infty \},
\]
\[
B^s_{\infty,q}(\R)=\{f\in C_b(\R): \|f\|^{(s)}_{\infty,q}<\infty \},
\]
\[
B^s_{p,\infty}(\R)=\{f\in L^p(\R):\omega_{k,p}(f;t)=O(t^s)\},
\]
\[
B^s_{\infty,\infty}(\R)=\{f\in C_b(\R):\omega_{k,\infty}(f;t)=O(t^s)\}.
\]
\end{definition}
\noindent One can prove that Besov spaces are independent of $k$.
If we would like to define  separable spaces we take $o(t^s)$ instead of $O(t^s)$.
\\
\\
Let us assume that a function $f$ belongs to $L^p(\R)$ space for some
$1\leq p<\infty$ or $f\in C_b(\R)$ in case $p=\infty$. We take the following definition of the smoothness parameter:
\begin{definition}\label{smooth}
 The value $s_p^*\in \mathbb{R}_+\cup\{0,\infty\}$ is the smoothness parameter of a function $f$ if
\begin{equation}
s_p^*=s_p^*(f)=\sup\{s>0: f\in B^s_{p,\infty}\},\nonumber
\end{equation}
where $\sup \emptyset =0$.
\end{definition}

\noindent By the following continuous embedding theorem
\[
B^{s_1}_{p,\infty}(\mathbb{R} )\subset B^{s_2}_{p,\infty}(\mathbb{R} ) \quad{\rm for}\quad s_1>s_2>0
\]
we obtain corollary:

\begin{corollary} The smoothness parameter $s^*_p$ of a function $f$ has the following properties
\begin{itemize}
\item $s^*_p=\infty$ $\Longleftrightarrow$ $f$ belongs to each $B^s_{p,\infty}(\mathbb{R} )$, where $s>0$

\item  $s^*_p=0$ $\Longleftrightarrow$ $f$ belongs to none of $B^s_{p,\infty}(\mathbb{R} )$, where $s>0$

\item  $0<s^*_p<\infty$  $\Longleftrightarrow$  for all $0<s<s^*_p$\ \ $f\in B^s_{p,\infty}(\mathbb{R} )$
and for all $s>s^*_p$\ \ $f\not\in B^s_{p,\infty}(\mathbb{R} )$.
\end{itemize}
\end{corollary}
\noindent The definition and facts about Besov spaces can be found in \cite{M} and \cite{HKPT}.
\\
\\
\noindent Below we give the characterization of $B^s_{p,\infty}(\mathbb{R} )$  in terms of a wavelet decomposition.
Following the notation in \cite{Hoffmann} let us define $r$-regular wavelet basis:

\begin{definition}\label{r-RWB}
Let $\phi$
be a scaling function and $\psi$ - the wavelet associated with $\phi$. If $r\geq 1$ is an integer then we will say that wavelet basis is $r$-regular if   $\phi , \psi \in C^r$ and the support of each of them is compact. We denote  $r$-regular wavelet basis by {\rm {$r$-RWB}}.
\end{definition}

\noindent There are many examples of $\phi$ and $\psi$
fulfilling the above definition (see \cite{D},  \cite{HKPT},\cite{M}, \cite{W}).
For instance, we can take Daubechies wavelets of a
sufficient large order.
Let us assume that
\begin{equation}\label{S(r)}
\supp \psi=[0,S(r)],
\end{equation}
where $S(r)\in \N$.
For a given $r$-RWB  we introduce some properties of $\psi$ and $\psi_{j,k}(\cdot):=2^{j/2}\psi(2^j\cdot -k)$.
Denote
\begin{equation}\label{PSI2}
\Psi_0:=\sup_{x\in \R} \sum_{k\in \Z}|\psi(x-k)|,\quad \Psi_2:=\sup_{x\in R} |\psi(x)| .
\end{equation}
By \cite[Corollary 5.5.2]{D} if we have $r$-regular wavelet basis then we have
the zero oscillation condition i.e. there is $d(r)\geq r$ such that
\begin{equation}\label{moment}
\int_{\mathbb{R}} x^k \psi(x) dx =0\quad {\rm for} \quad 0\leq k\leq d(r)
\quad \text{and} \ \
\int_{\mathbb{R}} x^{d(r)+1} \psi(x) dx \neq 0.
\end{equation}
Since we assume that $\psi$ is a real function then
\begin{equation}
\begin{aligned}
\label{periodic}
&\int_{0}^1 \left|\sum_{k\in\mathbb{Z}}\psi (u-k)\right|^2du  =\int_{0}^1 \sum_{k,j\in\mathbb{Z}}\psi (u-k)\psi (u-k-j) du\\
& =\sum_{k,j\in\mathbb{Z}}\int_{-k}^{-k+1}\psi (u)\psi (u-j) du =
\sum_{j\in\mathbb{Z}}\int_{\R} \psi (u)\psi (u-j) du=1.
\end{aligned}
\end{equation}
Let $I_{j,m}=[m/2^j,(m+1)/2^j]$. If $x\in I_{j,m}$, $y\in I_{j,l}$ then
\begin{equation}\label{kernel}
\sum_{k\in\mathbb{Z}}|\psi_{j,k}(x) \psi_{j,k}(y)|
\left\{
\begin{array}{lll}
= 0 & \text{for} & |m-l|>S(r)+1\\
\leq \Psi_0 \Psi_2 2^{j} & \text{for} & |m-l|\leq S(r)+1.
\end{array}
\right.
\\
\end{equation}
\noindent
Additionally with \eqref{S(r)}, we will need the following assumption to prove that our smoothness estimator is strongly consistent on the piecewise-smooth function class:
\begin{assumption}\label{assum1}
There exists $0<\delta_1<1$ such that $\psi (x)\neq 0$ for all $x\in (0,1+\delta_1]$.
\end{assumption}
\noindent This assumption is fulfilled for example by Daubechies wavelets DB2-DB20
(see  \cite[Lemma 3.1]{CD}). Note that since $\psi\in C^r$,  the assumption \ref{assum1} implies that there is $c>0$ such that
\begin{equation}\label{assum2}
\forall x\in [0,1] \quad \sum_{k\in\mathbb{Z}}\psi^2(x-k)>c> 0
\end{equation}
and
\begin{equation}\label{PSI}
\Psi_1:=\min_{0\leq n\leq d(r)} \left\{\left|\int_0^{\delta_1} (\delta_1-u)^n \psi(u) du\right| \right\}>0.
\end{equation}
By the assumption \ref{assum1} and the moment conditions \eqref{moment} we obtain that for all $\delta_1<\eta\leq 1+\delta_1$ and $0\leq n\leq d(r)$
\[
\left| \int_0^\eta (u-\eta)^n \psi(u) du \right|\geq \Psi_1.
\]
This property was crucial in a proof of \cite[Corollary 3.3]{CD}, and we will use it in this paper in Theorem \ref{p}.
\\
\\
Assuming that  $r$-RWB satisfies \eqref{assum2} we observe another property.
Since $\sum_{k\in\mathbb{Z}}\psi_{0,k}(x) \psi_{0,k}(y)$ is
 uniformly continuous then we conclude  by \eqref{assum2}
 that there are positive constants $a$ and $K$ such that if $|x-y|\leq 2^{-K}$ then
  $\left|\sum_{k\in\mathbb{Z}}\psi(x-k) \psi(y-k)\right|\geq a$. From this it follows that
\begin{equation}\label{iv}
\left|\sum_{k\in\mathbb{Z}}\psi_{j,k}(x) \psi_{j,k}(y)\right|\geq a2^j
\quad \text{for} \quad |x-y|\leq 2^{-(j+K)}.
\end{equation}
\\
\\
\noindent Now we introduce the wavelet decomposition of function $f$. For $j\geq 0$ we define a family of kernels
\begin{equation}\label{jadro}
G_j(x,y)=\sum_{k\in Z} \psi_{j,k}(x)\psi_{j,k}(y),
\end{equation}
and orthogonal projections
\begin{equation}\label{Qj}
Q_jf(x)=\int G_j(x,y) f(y) dy
\end{equation}
of function $f$. If we denote
\[
P_0f(x)=\sum_{k\in Z} <f,\phi(\cdot-k)> \phi(x-k)
\]
\begin{equation}\label{Pj}
P_{j+1}=Q_{j}+P_{j},\quad j\geq 0,
\end{equation}
where $<\cdot,\cdot>$ is the inner product in $L^2(\R)$, then for all
$f\in L^2(\R)$
\[
f=P_{0}f+\sum\limits_{j=0}^{\infty}Q_jf\ \ \ a.e.
\]
We have the following characterization of the Besov space  $B^s_{p,\infty}(\mathbb{R} )$ for $s<r$ and $1\leq p<\infty$ using $r$-RWB (see also for example  \cite{HKPT}, \cite{W})
\begin{equation}\label{char}
f\in B^s_{p,\infty}(\mathbb{R} ) \iff f\in L^p(\R ) \; {\rm and} \; \sup_{j\geq 0} 2^{js}\|Q_{j}f\|_p<\infty.
\end{equation}
For $p=\infty$ we have a similar characterization if we take $f\in C_b(\R )$.
Note that
\[
\|\psi_{jk}\|_p=2^{j/2-j/p}\| \psi\|_p.
\]
Since functions $\psi_{jk}, k\in \mathbb{Z}$ are orthonormal with a compact support we have the following stability condition ($1\leq p<\infty$): for all sequences $\{\beta_{jk}\}$
\begin{equation}\label{stability}
\frac{1}{\Psi_0} 2^{j/2-j/p}\left( \sum_{k\in Z} |\beta_{jk}|^p \right)^{1/p}\leq \left\|\sum_{k\in Z}\beta_{jk} \psi_{jk}\right\|_p \leq \Psi_0 2^{j/2-j/p} \left(\sum_{k\in Z} |\beta_{jk}|^p \right)^{1/p}.
\end{equation}
Let $\beta_{jk}(f)=\ <\psi_{jk},f>$. The following theorem was proved in \cite{CD} (see \cite[Theorem 2.1]{CD})
\begin{thm}
Let  {\rm $r$-RWB}  be given. Let $1\leq p\leq \infty$ and  $0<s^*_p<r$.  Then
\begin{equation}\nonumber
\liminf_{j\to \infty} \frac{-\log_2 \|\beta_{j\cdot}\|_{p}}{ j}=s^*_p-\frac{1}{p}+\frac{1}{2},
\end{equation}
where $\log_20=-\infty$.
\end{thm}
\noindent Applying the above theorem and \eqref{stability} we obtain corollary:
\begin{corollary}
\label{tw1} Let  {\rm $r$-RWB}  be given. Let $1\leq p\leq \infty$ and  $0<s^*_p<r$.  Then
\begin{equation}
\liminf_{j\to \infty} \frac{-\log_2 \|Q_j f\|_p}{ j}=s^*_p,\nonumber
\end{equation}
 where $\log_20=-\infty$.
\end{corollary}

\noindent In fact the above corollary is a simple consequence of the characterization \eqref{char}. We conclude that there is a subsequence $\{j_m\}=\{j_m(f)\}$
such that
\begin{equation}\nonumber
\lim_{m\to \infty} \frac{-\log_2^+ \|Q_{j_m} f\|_p}{ j_m}=s^*_p.
\end{equation}
In the next section we will see that if we assume that function $f$ is a piecewise-smooth function and a wavelet $\psi$ satisfies the assumption \ref{assum1} then we can replace $"\liminf"$ by $"\lim"$ in corollary \ref{tw1}. This is very important to obtain a strongly consistent estimator of the smoothness parameter and in a consequence a form of the smoothness test.
\\
\\

%%%%%%%%%%%%%%%%%%%%%%

\section{Piecewise-smooth functions}\label{psf}

In this section we focus on the piecewise-smooth functions class and the properties of the orthogonal projections $Q_j$ of those functions. Let us introduce the definition of the piecewise-smooth functions:

\begin{definition}\label{P-S}
A piecewise-smooth function $f$ with index $id(f)=m\geq 0$ is
a function with such properties:
\begin{itemize}
\item $f\in C^{m-1}(\R)$ if $m\geq 1$
\item There exist $n\in \mathbb{N}$ and $a_1,\ldots, a_n\in\mathbb{R}$ such that $f\notin C^{m}(\{a_i\})$ for some $i\in \{1,...,n\}$ \\and
$f\in C^{m+1}((-\infty,a_1])$, $f\in C^{m+1}([a_n,\infty))$,  $f\in C^{m+1}([a_i,a_{i+1}])$ for each i=1,...,n.
\end{itemize}
We say that the points $a_1,...,a_n$ are defect points of function $f$.
\end{definition}
\begin{remark}
An example of a piecewise-smooth function is any spline with finite number of knots. Notice that the orders of knots do not have to be the same. We say that the order of a knot $d\in\mathbb{R}$ of a spline $f$ is equal to $m$ if $f\in C^{m-1}(\{d\})$ for $m\geq 1$ and $f\notin C^{m}(\{d\})$.
The index of a spline function is the minimum of all orders of knots.
\end{remark}
\noindent Now we will prove several lemmas which we will need to prove the main theorem of this section. Recall that for a given {\rm $r$-RWB}, $d(r)$ is the number of vanishing moments of $\psi$.
\begin{lemma}\label{thm5}
 Let {\rm $r$-RWB} be given satisfying assumption \ref{assum1}.
Let $f$ be nonnegative   piecewise-smooth function with a compact
support and $id(f)\geq d(r)+1$. For $1<p<\infty$ we have
\begin{equation}\nonumber
\lim_{j\to \infty} 2^{pj(d(r)+1)}\int_{\mathbb{R}} |Q_jf(x)|^p f(x) dx
=\frac{1}{(d(r)+1)!^p}\int_0^1 |A(u)|^p du \int_{\mathbb{R}} |f^{(d(r)+1)}(x)|^p f(x) dx
\end{equation}
and
\begin{equation}\nonumber
\lim_{j\to \infty} 2^{pj(d(r)+1)}\int_{\mathbb{R}} |Q_jf(x)|^p  dx =\frac{1}{(d(r)+1)!^p}\int_0^1
|A(u)|^p du \int_{\mathbb{R}} |f^{(d(r)+1)}(x)|^p dx,
\end{equation}
where the function $A$ is defined by
\[
A(u):=\sum_{k\in \mathbb{Z}} \left( \int_{\mathbb{R}} (s-u)^{d(r)+1} \psi(s-k) ds \right) \psi(u-k).
\]
\end{lemma}
\begin{remark}
Notice that if $id(f)\geq d(r)+1$ then $f^{(d(r)+1)}$ and $f^{(d(r)+2)}$ exist almost everywhere so we can integrate those functions.
\end{remark}
\noindent \textbf{Proof}: The method of proof was developed in \cite{Dz} and \cite{BD}.
It is easily seen that $A$ is a periodic function, i.e. for all $l\in \mathbb{Z}$
\[
A(x)=A(x+l).
\]
First, let us assume that $f\in C^{d(r)+2}(\mathbb{R} )$ with a compact support. For each $x\in\mathbb{R}$ we can write Taylor's
polynomial
\[
\pi_x(u)=f(x)+f'(x)(u-x)+\cdots +\frac{f^{(d(r)+1)}(x)}{(d(r)+1)!}
(u-x)^{d(r)+1}
\]
and by the definition of $r$-RWB and \eqref{moment} we have
\begin{eqnarray}
Q_j(\pi_x)(x)&=&\sum_{k\in \mathbb{Z}} 2^j \left( \int_{\mathbb{R}}
\frac{f^{(d(r)+1)}(x)}{(d(r)+1)!} (u-x)^{d(r)+1} \psi(2^ju-k) du\right)
\psi(2^jx-k) \nonumber\\
&=& \frac{f^{(d(r)+1)}(x)}{(d(r)+1)!} {\frac{1}{ 2^{j(d(r)+1)}}} \sum_{k\in \mathbb{Z}}
\left(\int_{\mathbb{R}} (s-2^jx)^{d(r)+1} \psi(s-k) ds\right) \psi(2^jx-k)
\nonumber\\
&=&{\frac{1}{
2^{j(d(r)+1)}}} \frac{f^{(d(r)+1)}(x)}{(d(r)+1)!}A(2^j x).\nonumber
\end{eqnarray}

\noindent By Fejer-Orlicz-Mazur's theorem for periodic functions  (see \cite{BD}, \cite{Dz}) we obtain
\[
\lim_{j\to\infty}\int_{\mathbb{R}}
|f^{(d(r)+1)}(x)A(2^j x)|^p f(x) dx =\frac{1}{(d(r)+1)!^p} \int_{\mathbb{R}} |f^{(d(r)+1)}(x)|^p f(x)
dx \int_0^1 |A(u)|^p du,
\]
which gives
\[
\lim_{j\to\infty} 2^{pj(d(r)+1)} \int_{\mathbb{R}} |Q_j (\pi_x)(x)|^p f(x) dx =\frac{1}{(d(r)+1)!^p} \int_{\mathbb{R}} |f^{(d(r)+1)}(x)|^p f(x)
dx \int_0^1 |A(u)|^p du.
\]
\noindent By  \cite[Lemma 1.1]{Dz}
there exists $C>0$ such that for  all $f\in C^{d(r)+2}(\mathbb{R} )$ with
compact supports
\[
\int_{\mathbb{R}} |P_j (f-\pi_x)(x)|^p f(x) dx\leq C 2^{-pj(d(r)+2)} \int_{\mathbb{R}}
|f^{(d(r)+2)}(x)|^p f(x) dx,
\]
hence for some constant $C'$
\[
\int_{\mathbb{R}} |Q_j (f-\pi_x)(x)|^p f(x) dx\leq C' 2^{-pj(d(r)+2)} \int_{\mathbb{R}}
|f^{(d(r)+2)}(x)|^p f(x)dx,
\]
which completes the proof for $ f\in C^{d(r)+2}(\R)$.
Since $C^{d(r)+2}(\R)$ is dense in the Sobolev space $W^{d(r)+1}_p(\R)$, for all $1<p<\infty$ and any piecewise-smooth function with $id(f)\geq d(r)+1$ belongs to $W^{d(r)+1}_p(\R)$. The proof for the second formula is very similar.

\begin{lemma}\label{lem} If $f$ is nonnegative, piecewise-smooth function  with compact support and  $id(f)=k\in \mathbb{N}$ then there exists a constant  $C=C(f)>0$
such that
\[
\frac{1}{C} \int_{\mathbb{R}} |f^{(k)}(x)|^p  dx \leq \int_{\mathbb{R}} |f^{(k)}(x)|^p f(x)
dx \leq C \int_{\mathbb{R}} |f^{(k)}(x)|^p dx.
\]
\end{lemma}
\noindent \textbf{Proof}:  The right side of the above inequality is an easy
consequence of boundedness of the function $f$. To prove the left
side it is sufficient to make the following observation:
\[
\left( \int_{\mathbb{R}} |f^{(k)}(x)|^p f(x) dx=0 \right) \Rightarrow \left( \int_{\mathbb{R}} |f^{(k)}(x)|^p dx=0 \right)
\]
Suppose that
\[
\int_{\mathbb{R}} |f^{(k)}(x)|^p dx\neq 0.
\]
Since $f^{(k)}$ is piecewise-continuous there is an interval $(a,b)$ such that $|f^{(k)}(x)|>0$ for $x\in (a,b)$. But the function $f$ is nonnegative and continuous, so $f(x)>0$ for $x\in (a,b)$. Consequently,  $\int_{\mathbb{R}} |f^{(k)}(x)|^p f(x) dx>0$.
\\
\\
From the above lemmas we obtain the following corollary:
\begin{corollary}\label{cc}
Let $1<p<\infty$.
Let {\rm $r$-RWB} be given with $r\geq 1$. If $f$ is nonnegative,  piecewise-smooth function with compact support and  $id(f)\geq d(r)+1$, then there exist a natural number $N=N(f)$ and a constant $C>0$ such that for all
$j>N$
\begin{equation}
\label{Q}
\frac{1}{C}\int_{\mathbb{R}} |Q_jf(x)|^p dx \leq \int_{\mathbb{R}} |Q_jf(x)|^p f(x) dx \leq C
\int_{\mathbb{R}} |Q_jf(x)|^p  dx.
\end{equation}
\end{corollary}
\noindent \textbf{Proof}: Let the function $A(u)$ be defined as in lemma \ref{thm5}.
One can see that
\[
\int_0^1 |A(u)|^p du \neq 0.
\]
Indeed, from \eqref{moment} we obtain
\[
\int_{\mathbb{R}} (s-u)^{d(r)+1} \psi(s-k) ds  = \int_{\mathbb{R}} s^{d(r)+1}\psi (s)ds=: b\neq 0,
\]
hence
\[
A(u)=b\sum_{k\in\mathbb{Z}}\psi (u-k)
\]
and
\[
\int_0^1 |A(u)|^p du = b^p\int_0^1\left|\sum_{k\in\mathbb{Z}}\psi (u-k)\right|^p du >0
\]
according to \eqref{periodic}. Now using lemma \ref{lem} with $k=d(r)+1$ and lemma \ref{thm5} with $j$ large enough we have our result.
\\
\\
\begin{remark}
To obtain the main theorem of this section we will need one more theorem which was proved in \cite{CD} (see \cite[Corollary 3.3]{CD}). This theorem gives an asymptotic  characterization of $\|Q_j f\|_p$ for $1<p \leq \infty$. Moreover if we analyse the proof of the that Corollary    we will see that it is also true for $m\leq d(r)$, where $d(r)$ is the number of vanishing moments of $\psi$. In Proposition \ref{col}  we find the precise constants of estimates given in \cite[Corollary 3.3]{CD} for $1<p<\infty$. The case $p=\infty$ one can prove analogously.
\end{remark}

 Combining the above Remark  for $m\leq d(r)$ with Lemma \ref{thm5}, we obtain the full characterization of
 $\|Q_j f||_p$.  Moreover Corollary \ref{tw1} gives us relation between $id(f)$ and $s^*_p(f)$ for all $1<p\leq \infty$.
\begin{thm}
\label{p}
Let be given  {\rm $r$-RWB} satisfying assumption \ref{assum1}. Let $f$ be a piecewise-smooth function, bounded with a compact support. Then  we have
\begin{equation}\label{energia}
\|Q_j f\|_p\sim
\left\{
\begin{array}{lll}
2^{-j(id(f)+1/p)} & \text{if} & id(f)\leq d(r) \quad \text{and for}\quad 1< p\leq \infty\\
2^{-j(d(r)+1)} & \text{if} & id(f)\geq d(r)+1,\quad f\geq 0 \quad \text{and for} \quad 1< p<\infty\\
\end{array}
\right.,
\end{equation}
where by $\|Q_j f\|_p\sim 2^{-j \tau} $ we mean that there are $C_1,C_2$ and $N$ (dependent on $f$) such that for all $j\geq N$
\[
C_1  2^{-j\tau}\leq \|Q_j f\|_p\leq C_2   2^{-j\tau}.
\]
Moreover
\[
\lim_{j\to \infty} \frac{-\log_2^+ \|Q_j f\|_p}{ j}=
\left\{
\begin{array}{lll}
s^*_p(f) & \text{if} & id(f)\leq d(r)  \quad \text{and for} \quad 1< p\leq \infty\\
d(r)+1 & \text{if} & id(f)\geq d(r)+1, \quad f\geq 0 \quad \text{and for} \quad 1< p<\infty.\\
\end{array}
\right.
\]
\end{thm}

\begin{remark}
From the above theorem it is easy to see that if $f$ is piecewise-smooth function with a compact support then $s^*_p(f)=id(f)+1/p$ for all $1< p\leq  \infty$.
\end{remark} 	

\noindent To construct the form of the smoothness test we will need  more precise evaluation of $\|Q_j f\|_p$ than \eqref{energia}. For this purpose let us define the following class of functions:
\begin{definition}\label{class}
We say that function $f$ belongs to the class ${\mathcal S}([-1,1],d(r),\Delta_1,\Delta_2, N_S)$ if:
\begin{itemize}
\item $\supp f\subset [-1,1]$
\item $f$ is piecewise-smooth function with $id(f)=m\leq d(r)$
\item There exist  $\Delta_1>0$ and $\Delta_2>0$ such that for each defect point $d$ of function $f$
$\Delta_1\leq \Delta(f,d)\leq \Delta_2,$
where $\Delta(f,d):=|f^{(m)}(d^-)- f^{(m)}(d^+)|$
\item The number of defects is not greater than $N_S\in\mathbb{N}$
\end{itemize}
\end{definition}

\noindent If we assume that function $f$ belongs to the class ${\mathcal S}([-1,1],d(r),\Delta_1,\Delta_2, N_S)$ then we can obtain the following proposition:
\begin{Proposition}\label{col}
Let $1<p<\infty$. For all $f\in {\mathcal S}([-1,1],d(r),\Delta_1,\Delta_2, N_S)$ there exists $N=N(f)$ such that for all $j\geq N$
\[
\frac{1}{\Psi_0}   \frac{\Psi_1}{m!2^{j(m+1/p)}} \Delta_1-O(2^{-j(m+1)}) \leq \|Q_j f\|_p\leq   \Psi_0  \frac{\Psi_2 (S(r))^{m+2}}{(m+1)!2^{j(m+1/p)}}N_S \Delta_2 +O(2^{-j(m+1)}),
\]
where $m=id(f)$, constants $S(r),\Psi_0, \Psi_2, \Psi_1$ are defined in \eqref{S(r)},  \eqref{PSI2}, \eqref{PSI}.
\end{Proposition}
\noindent \textbf{Proof}:
This proposition can be proved in the same way as \cite[Corollary 3.3]{CD}. Let us define the function:
\[
h(x):=\frac{(-1)^{m}}{m!}D^-[(x_0-x)_+]^{m}+\frac{1}{m!}D^+[(x-x_0)_+]^{m},
\]
where constants $D^+,D^-$ are such that $D^+\neq D^-$.
Then for $m\leq d(r)$ we have
\begin{equation}\label{beta}
\frac{\Psi_1}{m!2^{j(m+1/2)}}|D^+-D^-| \leq\|\beta_{j\cdot}(h)\|_{l^p}\leq \frac{\Psi_2 (S(r))^{m+2}}{(m+1)!2^{j(m+1/2)}}|D^+-D^-|  ,
\end{equation}
where
\[
\|\beta_{j\cdot}(h)\|_{l^p}=\left(\sum_{k\in Z} |\beta_{j k}(h)|^p \right)^{1/p},\quad
\beta_{j k}(h)=\int\limits_{\mathbb{R}} h(x) \psi_{jk}(x)dx
\]
(see the proof of \cite[Corollary 3.3]{CD}).
Next for any $f\in f\in {\mathcal S}([-1,1],d(r),\Delta_1,\Delta_2, N_S)$    we can define functions
$$h_{m,i}(x):=\frac{(-1)^{m+1}}{m!}f^{(m)}(a_i^-)[(a_i-x)_+]^{m}-\frac{1}{m!}f^{(m)}(a_i^+)[(x-a_i)_+)]^{m},\ \ \ i=1,...,n$$
and
$$h_{m+1,i}(x):=\frac{(-1)^{m+2}}{m+1!}f^{(m+1)}(a_i^-)[(a_i-x)_+]^{m+1}-\frac{1}{m+1!}f^{(m+1)}(a_i^+)[(x-a_i)_+)]^{m+1},\ \ \ i=1,...,n$$
where the points $a_1,...,a_n$ are the defect points of function $f$. There exist function $H$ with a compact support, which equal to $h_{m,i}+h_{m+1,i}$ in some neighborhood of $a_i$ for each $i=1,...,n$ and $H\in C^{m+1}(\mathbb{R}/\{a_1,...,a_n\})$.
By adding $H$ to function $f$ we
remove all defects of function $f$ and obtain a function $g$ that belongs to Sobolev space $\ W^{m+1}_p(\R)$, for $1<p<\infty$ (see the proof of \cite[Corollary 3.3]{CD}).
By approximation properties of $P_j$ (see \cite[ Corollary 8.2]{HKPT}) we obtain
that for $g \in W^{m+1}_p(\R)$ and $m\leq d(r)$  there exists $C_{r,p}$ such that for all $j\geq 0$
\begin{equation}\label{sobolev}
\|Q_jg\|_p\leq C_{r,p} 2^{-j(m+1)}\|g^{(m+1)}\|_p.
\end{equation}
Note that by the triangle inequality
\[
\big|\|Q_jH\|_p-\|Q_j g\|_p\big|\leq \|Q_j f\|_p\leq \|Q_j g\|_p+\|Q_jH\|_p.
\]
But $\|Q_jH\|_p$ is controlled by $\|Q_j h_{m,i}\|_p$, since
$Q_j H(x)=Q_j h_{m,i}(x)+Q_j h_{m+1,i}(x) $
 in  a neighborhood of a point $a_i$, i.e.
 $x\in (a_i-\pi_i,a_i+\pi_i)$.
Consequently by  \eqref{stability}, \eqref{beta} and \eqref{sobolev} we obtain the result.

\begin{remark}
Note that for $p=2$ and  $f\in {\mathcal S}([-1,1],d(r),\Delta_1,\Delta_2, N_S)$ there exists $N=N(f)$ such that for all $j\geq N$
\begin{equation}\label{energia22}
 \frac{\Psi_1}{m!2^{j(m+1/2)}} \Delta_1 -O(2^{-j(m+1)})\leq \|Q_j f\|_2  \leq  \frac{\Psi_2 (S(r))^{d(r)+2}}{(m+1)!2^{j(m+1/2)}}N_S \Delta_2 +O(2^{-j(m+1)}),
\end{equation}
where $m=id(f)$, constants $S(r), \Psi_2, \Psi_1$ are defined in \eqref{S(r)},  \eqref{PSI2}, \eqref{PSI}.
\\
\end{remark}

%%%%%%%%%%%%%%%%%%

\section{Smoothness estimator}
Smoothness estimator was already introduced in \cite{Dziedziul} and \cite{CD}. In this section we refer to those results and present a slightly modified smoothness estimator which is more convenient for our purpose.
From now on, we will consider,  the case of $p=2$, i.e. a problem of estimating $s^*_2$.
As usual, let $X_1,X_2,...$ be a sequence of iid random variables with a density $f\in L^2(\mathbb{R} )$.
First, let us recall an estimator of $\|Q_j f\|_2^2$
\begin{equation}\label{ENJ}
{E_{n,j}}=\sum_{k\in \mathbb{Z}} \left|{\frac{1}
{n}}\sum_{i=1}^n \psi_{j,k}(X_i)
\right|^2
\end{equation}
which was examined in \cite{CD}.
One can see that
\[
{E_{n,j}}={\frac{2}
{n^2}} \sum_{i<l}^n G_j(X_i,X_l)+\frac{1}{n^2}\sum_{i=1}^n G_j(X_i,X_i).
\]
where $G_j$, for $j\geq 0$ are kernel functions defined in \eqref{jadro}.
Since for $i\neq k$
\begin{equation}\label{nieob}
\mathbf{E}[G_j(X_i,X_k)]=\sum_{k\in \mathbb{Z}} (\mathbf{E}\psi_{j,k}(X_1))^2= \sum_{k\in \mathbb{Z}}
(\langle\psi_{j,k},f\rangle )^2=\|Q_{j}f\|_2^2
\end{equation}
and
\[
\mathbf{E}[G_j(X_i,X_i)]=\sum_{k\in \mathbb{Z}} \mathbf{E}\psi^2_{j,k}(X_1)= \sum_{k\in \mathbb{Z}} \langle\psi^2_{j,k},f\rangle
\]
we have
\[
\mathbf{E}[E_{n,j}]=\frac{n-1}{n}\|Q_{j}f\|_2^2+\frac{1}{n} \sum_{k\in \mathbb{Z}} \langle\psi^2_{j,k},f\rangle
\]
so $E_{n,j}$ is biased estimator of $\|Q_{j}f\|_2^2$.
The following theorem was proved in \cite{CD} (see \cite[Corollary 4.2]{CD})
\begin{thm}\label{TwENJ}
Let be given  {\rm $r$-RWB} satisfying assumption \ref{assum1}.
Let $X_1,X_2,...X_n$ be a sequence of i.i.d random variables with density $f$, where $f$ is a piecewise-smooth function with $id(f)<r$, bounded and compactly supported. Then
\[
\lim_{n\to \infty} \frac{-\log_2 E_{n,j(n)}}{2j(n)}=id(f)+\frac{1}{2}=s^*_2\ \  a.e.,
\]
where $2^{j(n)(2r+1)}\asymp n$ (for example $j(n)=\lfloor\log_2n/(2r+1)\rfloor
$).
\end{thm}
\begin{remark}
The above theorem says that $-\log_2 E_{n,j(n)}/(2j(n))$ is a strongly consistent estimator of the smoothness parameter if $f$ is a piecewise-smooth function. If we analyze the proof of that theorem then, using Theorem \ref{p}, one can see that it is also true for $id(f)\leq d(r)$, where $d(r)$ is the number of vanishing moments of $\psi$.
\end{remark}
\noindent For the above construction of the smoothness estimator we have used a biased estimator of $\|Q_{j}f\|_2^2$. In the next section we will need an unbiased version of that estimator (U-estimator). We consider this estimator since it is easier to
formulate a concentration theorem for $L_{n,j}$. Let us introduce:
\begin{equation}\label{LNJ}
L_{n,j}={\frac{2}
{n(n-1)}} \sum_{i=1}^{n-1} \sum_{l=i+1}^{n} G_j(X_i,X_l).
\end{equation}
Using \eqref{nieob} we obtain that $L_{n,j}$ is an unbiased estimator of $\|Q_{j}f\|_2^2$. It is easy to see that
\[
L_{n,j}=\frac{n}{n-1}E_{n,j}-\frac{1}{n(n-1)}\sum_{i=1}^n G_j(X_i,X_i)
\]
By \eqref{kernel} there exists $C>0$ such that for all $j$
\[
G_j(x,x)=\sum_{k\in Z} \psi^2_{jk}(x) \leq C 2^j,
\]
then for  $2^{(2d(r)+1)j}\asymp n$ we have
\[
\left|\frac{1}{n(n-1)}\sum_{i=1}^n G_j(X_i,X_i)\right|\leq \frac{C 2^j}{n}\asymp 2^{-2jd(r))}\ \ a.e.
\]
Consequently we obtain the following corollary:
\begin{corollary}\label{colLNJ}
Let be given  {\rm $r$-RWB} satisfying assumption \ref{assum1}.
Let $X_1,X_2,...X_n$ be a sequence of i.i.d random variables with density $f$, where $f$ is a piecewise-smooth function with $id(f)\leq d(r)$, bounded and compactly supported. Then
\[
\lim_{n\to \infty} \frac{-\log_2 L_{n,j(n)}}{2j(n)}=id(f)+\frac{1}{2}=s^*_2\ \  a.e.,
\]
where $2^{j(n)(2d(r)+1)}\asymp n$ (for example $j(n)=\lfloor\log_2n/(2d(r)+1)\rfloor
$).
\\
\end{corollary}
\vspace{1.5cm}

%%%%%%%%%%%%%%%%%%%%%%%%%%%%%%%%%%%%%%%%%%%%

\section{Enrichment procedure}

In this section we focus on an additional condition for a density function $f$ which will be needed to construct a smoothness test in the next section. To avoid another restriction for the class of functions we introduce the enrichment procedure. First let us define our regularity condition:
\begin{definition}\label{regular} Let $f$ be  bounded with a compact support.
We say that a function $f$ is regular for given  $r$-RWB, if
 there are constants $N, C>0$ (dependent on $f$ and $r$-RWB) such that for all  $j=j(f)\geq N$
\begin{equation}
\label{ts1}
C \int_\R |Q_{j}f(x)|^2 dx \leq \int_{\mathbb{R}} (Q_{j} f(x))^2 f(x) dx.
\end{equation}
\end{definition}

\noindent In corollary \ref{cc} we have proved that if we have $r$-RWB then a piecewise-smooth function with $id(f)\geq d(r)+1$ satisfies condition \eqref{ts1}. On the other a hand a piecewise-smooth function with $id(f)\leq d(r)$ (which is in our area of interest) does not have to satisfy this condition. Let us introduce the following example:
\begin{example}\label{przyklad}
Let $r$-RWB be given with $r\geq 2$. Assume that
$\psi$ satisfies assumption \ref{assum1} and  ${\rm supp}\,\psi=[0,S(r)]$. Let us take the following density
\[
f(x)=\max\left( \frac{3}{4}(1-x^2),0\right).
\]
One can check that
\[
\|Q_j f\|_2 \approx 2^{-3j/2},
\]
i.e. $$s_2^*(f)=3/2.$$
From \eqref{moment} (zero oscillation condition) we have that for
$ |x-1|>S(r) 2^{-j}$ and $ |x+1|>S(r) 2^{-j}$ and every $j$ large enough
\[
Q_jf(x)=0.
\]
Since there is $C>0$ such that for all $j\geq 0$ and $x\in [-1;-1+S(r)2^{-j}]\cup  [1;1-S(r)2^{-j}]$
\[
f(x)\leq C2^{-j},
\]
then
\[
\int_{\mathbb{R}} (Q_j f(x))^2 f(x)dx \leq C 2^{-j} \int_{\mathbb{R}}(Q_j f(x))^2dx,
\]
so function $f$ is not regular in terms of definition \ref{regular}.
\end{example}

\noindent The above example shows that a problem appears when the smoothness of function is determined   only by defects $d_j$ such that $f(d_j)=0$. This leads us to the idea of the enrichment procedure.
\\
\\
Let $X_1,X_2,...,X_n$ be a sequence of iid random variables from a density  $f\in {\mathcal S}([-1,1],d(r),\Delta_1,\Delta_2, N_S)$.
If we add to that sample a new sequence of iid random variables $\tilde{X}_1,\tilde{X}_2,...,\tilde{X}_{n_1}$ from a density $\xi$, where
$
n_1=\pi n/(1-\pi),
$
and $0<\pi<1$, then we obtain a sample from the density
$
f_\pi=(1-\pi)f+\pi \xi,
$
of the size $n/(1-\pi)$.
Let
\begin{equation}\label{gfunkcja}
\xi(x)=\xi_\tau (x)=\frac{(2\tau+3)!}{[(\tau+1)!]^2}\ 3^{-(2r+3)}(1.5-x)^{\tau+1}(1.5+x)^{\tau+1}\ \mathbb{1}_{[-1.5;1.5]}(x).
\end{equation}
Note that  $\xi\in C^\tau(\R)$ and $\supp f \subset \supp \xi $. It is easy to see that for $\tau \geq d(r)+2$ function $f_\pi$ has the same smoothness as function $f$. Now we will show, by the following theorem, that function $f_\pi$ is regular in terms of definition \ref{regular}.

\begin{thm}\label{enrichment}
Let {\rm $r$-RWB} be given and $f\in {\mathcal S}([-1,1],d(r),\Delta_1,\Delta_2, N_S)$. Let $\xi=\xi_{d(r)+2}$ (see \eqref{gfunkcja}).
Then for  all  $0<\pi<1$ the function $f_\pi$ defined by formula
\[
f_{\pi}=(1-\pi) f +\pi \xi
\]
is regular with a constant $C=\pi\cdot\xi(1.25)/2,$
i.e. there exists $N=N(f,\xi,\pi)$ such that for $j\geq N$
\begin{equation}
\label{nienie1}
 \left(\int_{\mathbb{R}} (Q_{j} f_\pi(x))^2 f_\pi(x) dx\right)^{1/2} \geq
\sqrt{\pi\cdot\xi(1.25)/2}\ \|Q_{j}f_\pi\|_2.
\end{equation}
\end{thm}
\noindent \textbf{Proof}:
Let us fix $0<\pi<1$.
Recall
\[
\xi(x)=\xi_\tau (x)=\frac{(2\tau+3)!}{[(\tau+1)!]^2}\ 3^{-(2r+3)}(1.5-x)^{\tau+1}(1.5+x)^{\tau+1}\ \mathbb{1}_{[-1.5;1.5]}(x).
\]
There exists $\tilde{N}$ such that for all $j\geq \tilde{N}$ the length of ${\rm supp}\,\psi_{j,k}$ is smaller than $1/4$. Since  $\supp f\subset [-1,1]$ then for $j\geq \tilde{N}$ we have
\[
\pi^2 \int_{\mathbb{R}\setminus [-1.25,1.25]}|Q_j\xi(x)|^2dx =
\int_{\mathbb{R}\setminus [-1.25,1.25]}|Q_jf_{\pi}(x)|^2dx.
\]
Thus we have for $j$ large enough
\begin{eqnarray}
\lefteqn{\pi^2\int_{\mathbb{R}}|Q_j\xi (x)|^2dx +
\int_{\mathbb{R}}|Q_jf_{\pi}(x)|^2f_{\pi}(x)dx} \nonumber \\
& &\geq \pi^2\int_{\mathbb{R}\setminus [-1.25,1.25]}|Q_j\xi(x)|^2dx +
\int_{[-1.25,1.25]}|Q_jf_{\pi}(x)|^2f_{\pi}(x)dx \nonumber \\
& & \geq\int_{\mathbb{R}\setminus [-1.25,1.25]}|Q_jf_{\pi}(x)|^2dx +
\pi\ \xi(1.25)\int_{[-1.25,1.25]}|Q_jf_{\pi}(x)|^2dx \nonumber
\end{eqnarray}
Since $0<\xi(1.25)<1$ and $0<\pi<1$ then we get
\begin{equation}
\label{tsjm}
\pi^2 \|Q_j\xi\|^2_2 +
\int_{\mathbb{R}}|Q_jf_{\pi}(x)|^2f_{\pi}(x)dx
\geq  \pi\ \xi(1.25)\ \|Q_jf_{\pi}\|^2_2.
\end{equation}
Let $f\in {\mathcal S}([-1,1],d(r),\Delta_1,\Delta_2, N_S)$ and $\xi=\xi_{d(r)+2}$.
If $id(f)=m$ then $id(f_\pi)=m$. Hence
by \eqref{energia22}
\[
\|Q_j f_\pi \|_2 \sim 2^{-j(m+1/2)},
\]
and by \eqref{sobolev}
$$\|Q_j\xi\|_2\leq C_{r,2} 2^{-j(d(r)+1)} \|\xi^{(d(r)+1)}\|_2,$$
then  there exist $N$ such that for $j\geq N$
$$\pi^2\|Q_j\xi\|^2_2\leq \ \pi\ \xi(1.25)/2\ \|Q_j f_\pi \|^2_2.$$
Using \eqref{tsjm} and the above inequality we get our assertion.
\\
\\
\\
Now let us compare
the mean square error of the estimator with and without the enrichment procedure
in a numerical experiment. We use the function $f$
considered in the example \ref{przyklad}.
 The bias  is reduced significantly and variance
is smaller without enrichment procedure but we have a sample two times smaller!!
Note that $id(f)=1$ and $s^*_2(f)=3/2$.
We have generated $100$ samples each of the size $n=2^{20}$, thus we have calculated $100$ values of the smoothness estimator $-\log_2 L_{n,j(n)}/(2j(n))$ on the level
$j(n)=5$.

\begin{figure}[h!]
  \includegraphics[scale=0.5]{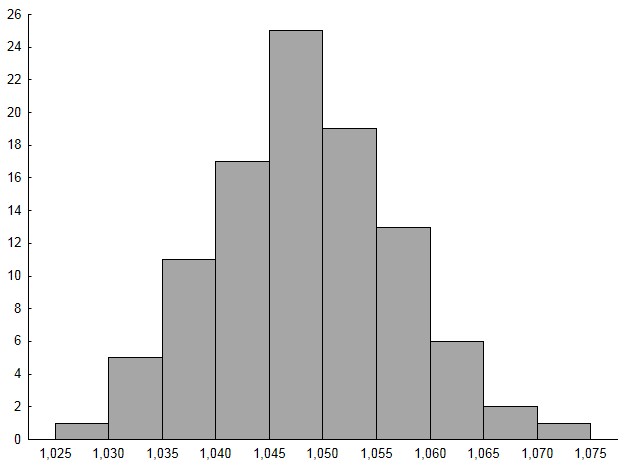}
\caption{Distribution of the smoothness estimator without enrichment procedure.}
\end{figure}

\noindent
Next we repeated our experiment using the procedure of sample enrichment with $\pi=\frac{1}{2}$ (for the results see Figure 2). It means that for each of 100 samples of size $ n=2^{20}$ we add a sample of size $n_1=2^{20}$ from the density $\xi_4$.

\begin{figure}[h!]
  \includegraphics[scale=0.5]{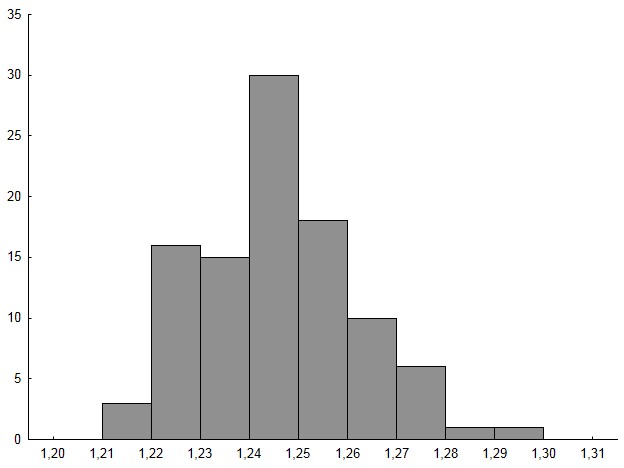}
\caption{Distribution of the smoothness estimator with enrichment procedure.}
\end{figure}

\noindent
It appears that the enrichment procedure gives better estimation of the smoothness parameter: in the second case, i.e. in the case of enrichment estimation, the mean value of the smoothness parameter was equal to
$1.25$  while in the first case $1.05$ when the true smoothness parameter is equal to $1.5$.
\\
\\
We also examine numerical results of changing $\pi$.
We enrich the old sample adding a sample from the density $\xi_4$ of the size $n_1=2^{20}\pi/(1-\pi) $. The influence of taking different values of $\pi$ on the smoothness estimator is presented on the figure below.

\begin{figure}[h!]
  \includegraphics[scale=0.5]{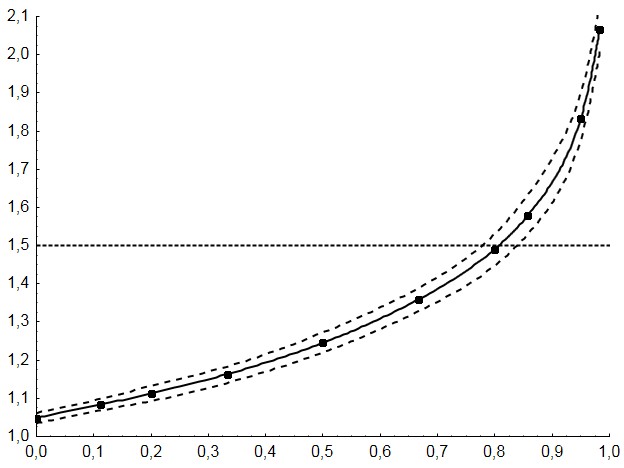}
\caption{Vertical axis: average of 100 smoothness estimators and approximate pointwise 95\% empirical confidence limits. Horizontal axis: $\pi$. The dashed horizontal line is the true value of the smoothness parameter.}
\end{figure}

\vspace{1.5cm}

%%%%%%%%%%%%%%%%%%%%

\section{Smoothness test}
To construct a smoothness test for a density function we will analyze the asymptotic distribution of the estimator $L_{n,j}$.
Let $X_1,X_2,...$ be a sequence of iid random variables with a density $f\in L^2(\mathbb{R} )$.
Denote
\[
\tilde{\sigma}_{j}:=\tilde{\sigma}_{j}(f)=\sqrt{\int_{\mathbb{R}} (Q_j f(x))^2 f(x) dx-\left(\int_{\mathbb{R}} (Q_j f(x))^2
dx\right)^2},
\]
\[
h_j(X_i,X_l):=G_j(X_i,X_l)-\|Q_{j}f\|_2^2,\ \ \ \ \ \textbf{Var}[h_j(X_1,X_2)]:=\sigma_j^2=\sigma_j^2(f)
\]
and
\[
g_j(x):=E[h_j(X_1,X_2)|X_1=x].
\]
 Then
\begin{eqnarray}
g_j(x)&=&E[G_j(X_1,X_2)|X_1=x]-\|Q_{j}f\|_2^2 =
\sum_{k\in Z}
\psi_{j,k}(x)E\psi_{j,k}(X_1)-\|Q_{j}f\|_2^2\nonumber \\
&=&Q_{j}f(x)-\|Q_{j}f\|_2^2.\nonumber
\end{eqnarray}

Note that the variance  of  $g_j(X_1)$
 is equal to  $\tilde{\sigma}_j^2=\tilde{\sigma}_j^2(f)$. Namely
\[
Eg_j^2(X_1)=E(Q_{j}f(X_1))^2-\|Q_{j}f\|_2^4
=\int_{\mathbb{R}} (Q_j f(x))^2 f(x) dx-\left(\int_{\mathbb{R}} (Q_j f(x))^2
dx\right)^2=\tilde{\sigma}^2_{j}.
 \]
To find an asymptotic formula for $\sigma_j$ we introduce the following lemma:
\begin{lemma}\label{s}
Let be given  {\rm $r$-RWB} satisfying assumption \ref{assum1}.
Let $X_1,X_2,...$ be a sequence of iid random variables with a density $f\in L^2(\mathbb{R} )$. Then there are constants $N, C_1, C_2>0$ independent of $f$ such that for each $j\geq N$
\begin{equation}\nonumber
C_1 \|f\|_2^2\leq2^{-j}\sigma_j^2\leq C_2 \|f\|_2^2.
\end{equation}
\end{lemma}
\noindent{\bf Proof}:
Note that
\[
\sigma_j^2(f)=E\left( G_j(X_1,X_2)\right) ^2-\|Q_jf\|_2^4.
\]
Since $\|Q_jf\|_2^4\to 0$, we only need to show an appropriate estimate on the term $E\left( G_j(X_1,X_2)\right) ^2$. Let $I_{j,m}=[m/2^j,(m+1)/2^j]$.
By \eqref{kernel} we obtain
\begin{eqnarray}
E\left( G_j(X_1,X_2)\right) ^2 &=& \int_{\mathbb{R}}\int_{\mathbb{R}}\left(\sum_{k\in\mathbb{Z}}\psi_{j,k}(x)\psi_{j,k}(y)\right) ^2f(x)f(y)dxdy \nonumber \\
&\leq& C\sum_{m\in\mathbb{Z}}\sum_{l:\; |m-l|\leq D}\int_{I_{j,m}}\int_{I_{j,l}}2^{2j}f(x)f(y)dxdy \nonumber \\
&=&C2^{2j}\sum_{m\in\mathbb{Z}}\left( \int_{I_{j,m}}f(x)dx\right)\left( \sum_{l:\; |m-l|\leq D}\int_{I_{j,l}}f(y)dy\right)\nonumber \\
&\leq& C'2^{2j}\sum_{m\in\mathbb{Z}}\left( \int_{I_{j,m}}f(x)dx\right)^2\leq C'2^{j}\|f\|_2^2.\nonumber
\end{eqnarray}
On the other hand, by \eqref{iv} we conclude that
\begin{eqnarray}
E\left( G_j(X_1,X_2)\right) ^2 &=& \int_{\mathbb{R}}\int_{\mathbb{R}}\left(\sum_{k\in\mathbb{Z}}\psi_{j,k}(x)\psi_{j,k}(y)\right) ^2f(x)f(y)dxdy\nonumber \\
&\geq& \sum_{m\in\mathbb{Z}}\int_{I_{j+K,m}}\int_{I_{j+K,m}}\left(\sum_{k\in\mathbb{Z}}\psi_{j,k}(x)\psi_{j,k}(y)\right) ^2f(x)f(y)dxdy\nonumber \\
&\geq& \sum_{m\in\mathbb{Z}}\int_{I_{j+K,m}}\int_{I_{j+K,m}}2^{2j}a ^2f(x)f(y)dxdy \nonumber \\
&=& 2^{2j}a^2\sum_{m\in\mathbb{Z}}\left( \int_{I_{j+K,m}}f(x)dx\right)^2\geq
C''2^j\|f\|_2^2, \nonumber
\end{eqnarray}
which is our claim.
\\
\\
To find an asymptotic distribution of the estimator $L_{n,j}$ we will also need the following lemma:
\begin{lemma}\label{end}
Let $r$-RWB be given. Let $f$ be a bounded density and $\supp f\subset [-1,1]$.
Then there are constants $C=C(\|f\|_\infty)$ and $C'=C'(\|f\|_\infty)$ such that for all $j\geq 0$
\[
E|g_j(X_1)|^3\leq C \int_\mathbb{R} |Q_jf(x)|^3dx\leq C' \|Q_j f\|_2^2.
\]
\end{lemma}

\noindent {\bf Proof}:
Since $\supp f\subset [-1,1]$ and $\supp \psi=[0,S(r)]$ then there is an interval $[a,b]$ such that for all $j\geq 0$
\begin{equation}\label{one}
 {\rm supp}\; Q_jf\subset [a,b].
\end{equation}
Hence
\begin{equation}\label{two}
| Q_jf(x)|\leq  \sum_{k}\langle |\psi_{j,k}|,f \rangle |\psi_{j,k}(x)| \leq \|f\|_{\infty}\cdot\|\psi\|_{1}\cdot\|\psi\|_{\infty}\cdot (b-a).
\end{equation}
Let
\[
\Psi_4:=\|\psi\|_{1} \|\psi\|_{\infty}  (b-a).
\]
For $1\leq p <\infty$ by \eqref{one} and \eqref{two} we have that for all $j\in\N$
\[
\| Q_jf\|_{p}=\left( \int_{\R} |Q_jf(x)|^p\, dx\right)^{1/p}=\left( \int_a^b |Q_jf(x)|^p\, dx\right)^{1/p}.
\]
Consequently
\begin{equation}\label{three}
\| Q_jf\|_{p}\leq \left( \int_a^b (\Psi_4 \|f\|_\infty)^p\, dx\right)^{1/p}=\Psi_4 \|f\|_\infty (b-a)^{1/p}.
\end{equation}
For each $1\leq p_1 < p_2 <\infty$ by  H\"{o}lder inequality we get
\begin{eqnarray}
\| Q_jf\|_{p_1}^{p_1}&=& \int_a^b |Q_jf(x)|^{p_1}\, dx
\leq \left( \int_a^b \left(|Q_jf(x)|^{p_1}\right)^{p_2/p_1}\, dx\right)^{p_1/p_2}\cdot \left(\int_a^b1\, dx\right)^{1-p_1/p_2}\nonumber \\
&=& (b-a)^{1-p_1/p_2}\left(\| Q_jf\|_{p_2}\right)^{p_1}, \nonumber
\end{eqnarray}
which gives
\begin{equation}\label{four}
\| Q_jf\|_{p_1}\leq (b-a)^{\frac{p_2-p_1}{p_1p_2}}\| Q_jf\|_{p_2}.
\end{equation}
Now we are ready to evaluate $E|g_j(X_1)|^3$.
From the obvious inequality $|s-t|^3\leq 7(|s|^3+|t|^3)$ and the fact that $f$ is a bounded density we get
\begin{eqnarray}
E|g_j(X_1)|^3&=&\int_{\R}\left| Q_jf(x)-\|Q_jf\|_2^2\right|^3f(x)\, dx \nonumber \\
&\leq& 7\left(\int_{\R}| Q_jf(x)|^3f(x)\, dx + \|Q_jf\|_2^6\int_{\R}f(x)\, dx\right) \nonumber \\
&\leq& 7\max\{\|f\|_\infty,1\}\left(\int_{\R}| Q_jf(x)|^3\, dx + \|Q_jf\|_2^6\right). \nonumber
\end{eqnarray}
Using \eqref{three} and \eqref{four} we have
\[
\|Q_jf\|_2^6\leq (b-a)\|Q_jf\|_3^6=(b-a)\|Q_jf\|_3^3 \|Q_jf\|_3^3
\leq (b-a)\|Q_jf\|_3^3 (\Psi_4 \|f\|_\infty)^3 (b-a).
\]
Finally, there is $C=C(\|f\|_\infty)>0$ such that all $j\geq 0$
\[
E|g_j(X_1)|^3\leq C  \cdot \int_{\R}| Q_jf(x)|^3\, dx.
\]
Moreover, by \eqref{two} the following evaluation is true
\[
\int_{\R}| Q_jf(x)|^3\, dx\leq \Psi_4 \|f\|_\infty \int_{\R}| Q_jf(x)|^2\, dx,
\]
which finishes the proof.
\\
\\
The final fact we will need is a classical theorem for U-statistics. Let $\Phi$ denote the  standard  Normal distribution function.

\begin{thm}(Berry Esseen inequality)  Let $X_1,X_2,\ldots , X_n$ be a sequence of iid random
variables and let $U_n$ be given by
\[
U_n={\frac{2}
{n(n-1)}} \sum_{1\leq l<s\leq n} h(X_l,X_s),
\]
where $h(x,y)$ is a symmetric, real-valued function. Let $Eh(X_1,X_2)=0$, $\sigma^2=Eh^2(X_1,X_2)<\infty$ and $\tilde{\sigma}^2=Eg^2(X_1)>0$ where $g(x)= E(h(X_1,X_2)|X_1=x)$. If $(E|g(X_1)|^3<\infty$ then
\[
\sup_{z\in\mathbb{R}}\left| P\left(\frac{\sqrt{n}}{ 2 \tilde{\sigma}}U_{n}\leq
z\right)-\Phi(z)\right|\leq \frac{6.1 E|g(X_1)|^3}
{\sqrt{n}\tilde{\sigma}^3}+\frac{(1+\sqrt{2})\sigma}
{\sqrt{2(n-1)}\tilde{\sigma}}.
\]
\end{thm}
\noindent The proof of this theorem, as well as much more details on U-statistics, can be found for example in \cite{CGS}.
\\
\\
Now we can formulate our main theorem which allow us to construct a smoothness test for a density function that belongs to the class ${\mathcal S}([-1,1],d(r),\Delta_1,\Delta_2, N_S)$:
\begin{thm}\label{thm55}
Let $r$-RWB be given satisfying assumption \ref{assum1}.
Let $X_1,X_2,...,X_n$ be a sequence of iid random
variables with density $f\in {\mathcal S}([-1,1],d(r),\Delta_1,\Delta_2, N_S)$.  Let us use the enrichment procedure for $\xi_{d(r)+2}$ and fixed $0<\pi<1$.  There is $C>0$ such that for all $f_\pi$,  there is a constant  $N\in\mathbb{N}$ such that for all $j\geq N$
\begin{equation}\nonumber
\sup_{z\in\mathbb{R}}\left( P\left({\frac{\sqrt{n}}{ 2 \tilde{\sigma}_{j}(f_\pi)}}\left[ L_{n,j}-\|Q_{j}f_\pi\|_2^2\right]\leq
z\right)-\Phi(z)\right)^2\leq C n^{-\frac{1}{2m+3}},
\end{equation}
where
\[
n\geq C_1 2^{j(2m+3)},
\]
and $id(f)=m$, $C_1>0$.
Note that $n$ is a new sample size (after enrichment).
\end{thm}
\noindent\textbf{Proof}:
By lemma \ref{end}  we have
\[
E|g_j(X_1)|^3\leq C_2 \int_R |Q_jf(x)|^3dx\leq C_3 \|Q_j f\|_2^2.
\]
Let us fix $\pi$. Now  we can use \eqref{nienie1}. Moreover there is $\Psi_5>0$ such that for all $f\in {\mathcal S}([-1,1],d(r),\Delta_1,\Delta_2, N_S)$ we have $\|f\|_\infty \leq \Psi_5$. Consequently the constant in the evaluation in lemma \ref{end} may be chosen independently from $f_\pi$. By Theorem \ref{enrichment} and Lemma \ref{s} we obtain that there is $C_4>0$  for $f_\pi$ such that if $j$ is large enough
\begin{equation}\label{berry}
\frac{6.1 E|g_{j}(X_1)|^3}{
n^{1/2}\tilde{\sigma}_{j}(f_\pi)^3}+\frac{(1+\sqrt{2})\sigma_{j}}{
\sqrt{2(n-1)}\tilde{\sigma}_{j}(f_\pi)}
\leq
\frac{C_4\sqrt{2^{j}}}{
\sqrt{n}\|Q_{j} f_\pi\|_2}.
\end{equation}
Using Proposition \ref{col} we obtain
\begin{equation}\label{energia25}
(1-\pi) \frac{\Psi_1}{m!2^{j(m+1/2)}} \Delta_1 - O(2^{-j(m+1)})\leq \|Q_j f_{\pi}\|_2  \leq  (1-\pi)\frac{\Psi_2 (S(r))^{d(r)+2}}{(m+1)!2^{j(m+1/2)}}N_S \Delta_2+ O(2^{-j(m+1)})
\end{equation}
for sufficient large $j\geq N$. By Berry Esseen inequality, \eqref{berry} and \eqref{energia25} we finish the proof.
\\
\\
\\
Now we are ready to construct the form of the following smoothness test for a density function $f\in{\mathcal S}([-1,1],d(r),\Delta_1,\Delta_2, N_S)$  for fixed  $\mu\in\{0,1...,d(r)\}$ :
\begin{equation}\nonumber
  H_0: s^*_2(f)\leq  \mu+1/2,\quad \text{against} \quad H_1: s^*_2(f)>  \mu+1/2.
\end{equation}
\begin{remark}
It is easy to see that if $f$ belongs to ${\mathcal S}([-1,1],d(r),\Delta_1,\Delta_2, N_S)$ then its smoothness parameter $s^*_2$ is equal to $m+1/2$ where $m=id(f)$ belongs to $\{0,1...,d(r)\}$.
\end{remark}
\noindent Let $X_1,X_2,...,X_n$ be a sequence of iid random
variables with density $f$.
Since we want to use Theorem \ref{thm55} we enrich our sample for a given $0<\pi<1$. If we denote by $z_\alpha$ the quantile of the standard normal distribution then we reject the null hypothesis for $n\asymp 2^{j(n)(2d(r)+3)}$ at the significance level $\alpha$ when
\begin{equation}\label{test}
L_{n,j(n)}\leq 2 \tilde{\sigma}_{j(n)}(f_{\pi})z_\alpha/\sqrt{n}+\|Q_{j(n)}f_{\pi}\|_2^2.
\end{equation}
By the definition of $\tilde{\sigma}_j$
\[
\sqrt{ \int_{\mathbb{R}} (Q_{j} f_{\pi}(x))^2 f_{\pi}(x) dx} \sim
 \tilde{\sigma}_j(f_{\pi}).
\]
Now by \eqref{nienie1} there exists $N$ such that for $j\geq N$
\begin{eqnarray}\label{test1}
\sqrt{\pi \xi_{d(r)+2}(1.25)/2}\ \|Q_j f_{\pi}\|_2 \leq \tilde{\sigma}_j(f_{\pi}).
\end{eqnarray}
Using \eqref{test}, \eqref{test1} and \eqref{energia25} we obtain the principle which we can use in practice.\\
\\
 \textbf{Smoothness test:}\\
We reject the null hypothesis for $ f\in{\mathcal S}([-1,1],d(r),\Delta_1,\Delta_2, N_S)$ at the significance level $\alpha$ when
\begin{eqnarray}\label{test2}
L_{n,j(n)} \leq
z_\alpha \frac{ \Delta_1  \Psi_1  \sqrt{\pi \xi_{d(r)+2}(1.25)/2}\ (1-\pi)}{\mu! \sqrt{n}\ 2^{j(n)(\mu+1/2)}}+\left(\Delta_1  \frac{\Psi_1}{\mu!}\right)^2 (1-\pi)^2 2^{-j(n)(2\mu+1)},
\end{eqnarray}
where $\Psi_1$, $\xi_{d(r)+2}$ are defined in \eqref{PSI} and \eqref{gfunkcja}, $j(n)=\lfloor\log_2n/(2d(r)+3)\rfloor$, $d(r)$ is the number of the vanishing moments of the wavelet and
$\pi$ is a parameter of the enrichment procedure.

\section{Numerical experiment}
In the proofs of the above theorems the exact values of constants are not needed since we examine asymptotic formulas.
In applications  we  need exact constants. One can calculate them numerically. Unfortunately one of our constants is not very good in applications. The constant $\Psi_1$ (see (\ref{PSI})) is very small for Daubechies wavelets  which affects on a rejection area of our test for small resolution level $j(n)$. Since constant $\Psi_1$ is very comfortable in proofs but not very useful in applications we suggest to change it depending on the type of the test.
Let us consider the following example:
\begin{example}
 Let $\psi=DBN$ be a Daubechies wavelet ($DBN$ where $N\in\{1,2,...,20\}$). Let  $h$ be a characteristic function of an interval $I$, i.e. $h=\chi_I$ such that $|I|=1>2^{-j_0}\supp \psi$. Since we have two points of discontinuity then  for all $j\geq j_0$ (see the  proof of Lemma 3.2 and Remark 3.2 \cite{CD})
$$
\frac{\log_2 2 \sup F_\psi}{2j}\geq  -\frac{\log_2 \|Q_j h\|_2^2}{2j}-\frac{1}{2}- id(h) \geq \frac{\log_2 2 \inf F_\psi}{2j},
$$
where the function $F_\psi$ is given numerically (see Figure 3, \cite{CD}). For instance, let $\psi=DB8$. Then $\sup F_\psi\approx 0.08$ and $\inf F_\psi\approx 0.02$.  So
$$
 \left| -\frac{\log_2 \|Q_j h\|_2^2}{2j}-\frac{1}{2}- id(h)\right|\leq -\frac{\log_2 0.02 }{2j}\approx \frac{5,64}{2j} .
$$
We can see that detection of discontinuity of function $h$ in non random case require at least $j=5$. Since $n\asymp 2^{j(2d(r)+3)}$ in random case, the sample size should be huge.
\end{example}
\noindent The example above shows that we can change $\Psi_1$ to $\sqrt{\inf F_\psi}$ if we want to test  $H_0: id(f)\leq  0$ against $H_1: id(f)\geq  1$. Furthermore, instead of $\Psi_1$ we can take a value of some sequence which converges to $\sqrt{\inf F_\psi}$ when $j \rightarrow \infty$. We suggest the following correction of the constant $\Psi_1$: $V_{j(n)}:= \sqrt{\inf F_\psi + 1/j(n)}$. Now let us check the behaviour of our test in the following numerical experiment:
\\
\\
Let us consider the following density functions:
$$f_0(x)=\mathbb{1}_{[0,1]}(x)\left(\frac{1}{2}+3x(1-x)\right),\ \ \ \ f_1(x)=\mathbb{1}_{[0,1]}(x)6x(1-x)$$

\begin{figure}[h!]
  \includegraphics[scale=0.4]{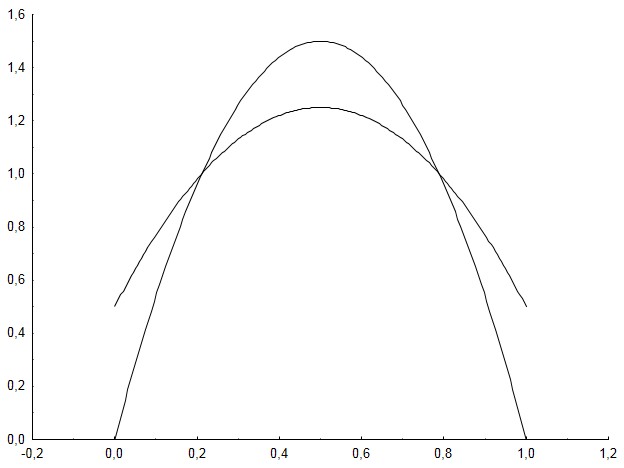}
\caption{Densities $f_0$ and $f_1$}
\end{figure}
\noindent  It is easy to see that $id(f_0)=0$ and $id(f_1)=1$. Three values of the experiment size were used: $n_1=2^{16}$, $n_2=2^{20}$, and $n_3=2^{24}$. The samples were enriched using function $\xi_3$ (see (\ref{gfunkcja})) with $\pi=1/2$. For the estimation the Daubechies wavelet DB8 was used (with support length $15$). The resolution levels were: $j(n_1)=4$, $j(n_2)=5$ and $j(n_3)=6$. Using (\ref{test2}) and the correction $V_{j(n)}$ the following rejection procedure was taken: We reject $H_0: id(f)\leq  \mu_0=0$ if
$$
\widehat{id(f)}_{n}=-\frac{\log_2 L_{n,j}}{2j} - \frac{1}{2}\geq \mu_0 -\frac{\log_2\left(\frac{z_\alpha \Delta_1  V_j  \sqrt{\pi \xi_{3}(1.25)/2}\ (1-\pi) 2^{j (\mu_0+1/2)}}{\mu_0! \sqrt{n}}+\left(  \frac{\Delta_1V_j(1-\pi)}{\mu_0!}\right)^2  \right)}{2j},
$$
where $\alpha=0.05$, $z_{\alpha}=-1.65$, $\Delta_1=1/2$, $ V_j=\sqrt{0,02 + 1/j}$, $\pi=1/2$, $\xi_{3}(1.25)=0,007$ and $\mu_0!=1$. The results are presented in the figures below:
\\
\\

\begin{figure}[h!]
  \includegraphics[scale=0.5]{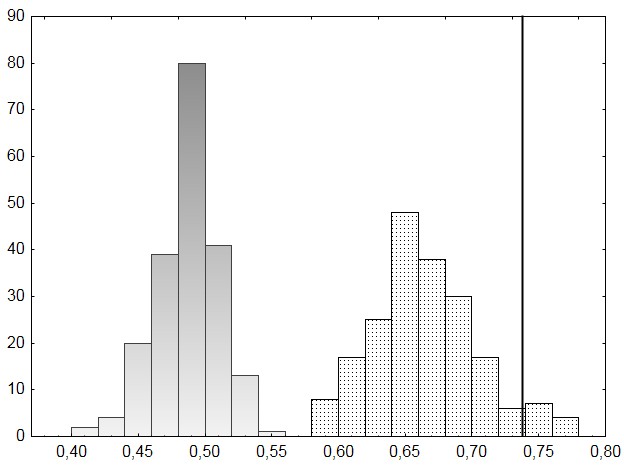}
\caption{Grey histogram presents 200 estimators of the index of function $f_0$ (enriched by function $\xi_3$ with $\pi=1/2$) from 200 generations with sample size $n_1=2^{16}$, and resolution level $j(n_1)=4$. Histogram filled with dots presents the same but for function $f_1$. The black vertical line determines the rejection area which is on the right side of that line.}\label{h200}
\end{figure}
\vspace{0.5cm}
\noindent One can see that in this example the test does not reject the null hypothesis when it is true, but also does not reject it in the most cases when it is false. It means that the power of our test for the sample size $n_1=2^{16}$ is very low. Now let us take the sample size $n_2=2^{20}$ and resolution level $j(n_2)=5$

\begin{figure}[h!]
  \includegraphics[scale=0.5]{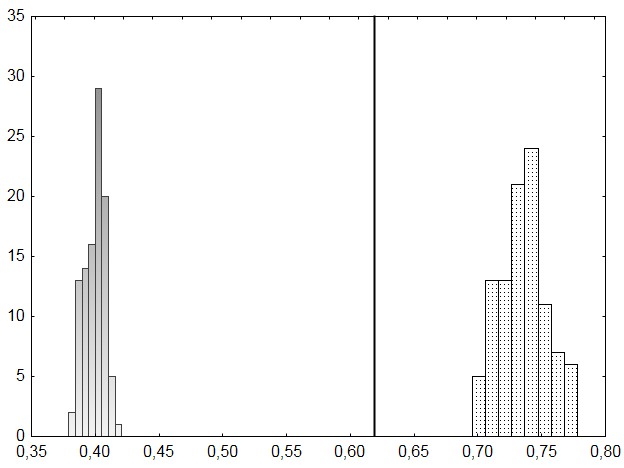}
\caption{Similar to Figure \ref{h200} but here we have 100 estimators from 100 generations with sample size $n_2=2^{20}$, and resolution level $j(n_2)=5$.}
\end{figure}
\noindent Now one can see that in all cases our test rejects the null hypothesis when it is false and does not reject it when it is true. It means that the empirical power for the sample size $n_1=2^{20}$, and resolution level $j(n_1)=5$ is equal to $1$. One can also see that the variance as well as the bias of the index estimator are smaller than in the previous case. Now let us check
what happens for the sample size $n_3=2^{24}$ and resolution level $j(n_3)=6$

\begin{figure}[h!]
  \includegraphics[scale=0.5]{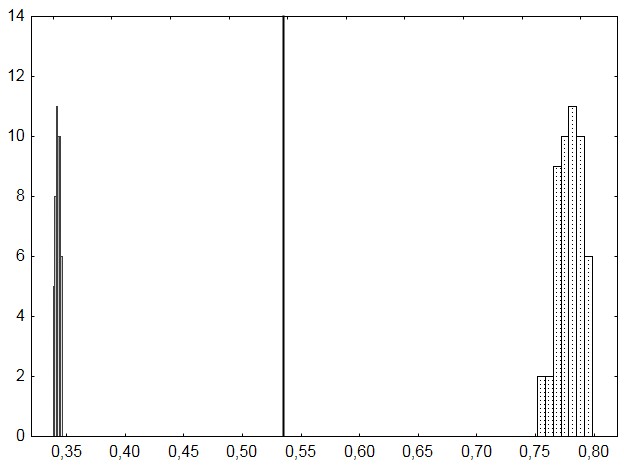}
\caption{Similar to Figure \ref{h200} but here we have 50 from 50 generations with sample size $n_3=2^{24}$, and resolution level $j(n_3)=6$.}
\end{figure}
\noindent One can see that the variance and the bias of the index estimator are even smaller than in the previous case. As in the previous case we do not observe the Type I and Type II errors.
\\
\\
\textbf{Acknowledgments}\\
The authors would like to thank the anonymous reviewers for their valuable comments and suggestions to improve the paper.

%%%%%%%%%%%%%%%%%%%%%%%%%

%%%%%%%%%%%%%%%%%%%%%%%%%

\newpage

\noindent Bogdan \'Cmiel\\
Faculty of Applied Mathematics\\
AGH University of Science and Technology\\
Al. Mickiewicza 30\\
30-059 Cracow\\
Poland\\
cmielbog@gmail.com \\
\\
\\
Karol Dziedziul\\
Faculty of Applied Mathematics \\
Gda\'{n}sk University of Technology\\
ul. G. Narutowicza 11/12\\
80-952 Gda\'nsk \\
Poland\\
kdz@mifgate.pg.gda.pl\\
\\
\\
Barbara Wolnik\\
Institute of Mathematics\\
University of Gda\'nsk\\
ul. Wita Stwosza 57\\
80-952 Gda\'nsk\\
Poland\\
Barbara.Wolnik@mat.ug.edu.pl

\end{document}